\documentclass[11pt,a4paper]{amsart}
\usepackage{graphicx,multirow,array,amsmath,amssymb}
\usepackage{xcolor}
\usepackage{longtable}
\usepackage{supertabular} 

\newtheorem{theorem}{Theorem}
\newtheorem{proposition}[theorem]{Proposition}
\newtheorem{lemma}[theorem]{Lemma}
\newtheorem{corollary}[theorem]{Corollary}
\theoremstyle{definition}
\newtheorem*{definition}{Definition}
\theoremstyle{remark}
\newtheorem*{remark}{Remark}

\begin{document}
	
\title{Linear Upper Bounds on the Ribbonlength of Knots and Links}
	
\author[H. Kim]{Hyoungjun Kim}
\address{Institute of Statistics, Korea University, Seoul 02841, Korea}
\email{kimhjun@korea.ac.kr}

\author[S. No]{Sungjong No}
\address{Department of Mathematics, Kyonggi University, Suwon 16227, Korea}
\email{sungjongno@kgu.ac.kr}

\author[H. Yoo]{Hyungkee Yoo}
\address{Department of Mathematics Education, Sunchon National University, Sunchon 57922, Korea}
\email{hyungkee@scnu.ac.kr}

\keywords{ribbonlength, folded ribbon knot}
\subjclass[2020]{57K10}
\thanks{
The first author(Hyoungjun Kim) was supported by the National Research Foundation of Korea (NRF) grant funded by the Korea government Ministry of Science and ICT(NRF-2021R1C1C1012299 and NRF-2022M3J6A1063595).
}
\thanks{
The second author(Sungjong No) was supported by the National Research Foundation of Korea(NRF) grant funded by the Korea government Ministry of Science and ICT(NRF-2020R1G1A1A01101724).
}
\thanks{
The corresponding author (Hyungkee Yoo) was supported by Basic Science Research Program of the National Research Foundation of Korea (NRF) grant funded by the Korea government Ministry of Education (RS-2023-00244488).
}

\begin{abstract}
A knotted ribbon is one of physical aspect of a knot.
A folded ribbon knot is a depiction of a knot obtained by folding a long and thin rectangular strip to become flat.
The ribbonlength of a knot type can be defined as the minimum length required to tie the given knot type as a folded ribbon knot.
The ribbonlength has been conjectured to grow linearly or sub-linearly with respect to a minimal crossing number.
Several knot types provide evidence that this conjecture is true, but there is no proof for general cases.
In this paper, we show that for any knot or link, the ribbonlength is bounded by a linear function of the crossing number.
In more detail,
$$
\text{Rib}(K) \leq 2.5 c(K)+1.
$$
for a knot or link $K$.
Our approach involves binary grid diagrams and bisected vertex leveling techniques.
\end{abstract}
	
\maketitle

\section{Introduction} \label{sec:intro}

Knot theory, a branch of topology, studies knots and their properties.
A {\it knot\/} is a simple closed curve in three dimensional space, and a {\it link\/} is a disjoint union of knots.
Understanding the properties of knots is useful in the wide area of natural science such as physics, biology, chemistry and materials science.
One physical aspect of knot theory is the concept of a ribbon structure of a knot.
A knotted ribbon represents a geometric measure of its complexity, providing valuable insights into its topological structure~\cite{CDBG,RLRA}.
A ribbon structure realizes a ribosomal walking robot~\cite{RCL,WPF}.
Many circular DNA~\cite{LJ,LLA,LPCW,SCBLC} have knotted ribbon shapes.

\begin{figure}[h!]
	\includegraphics{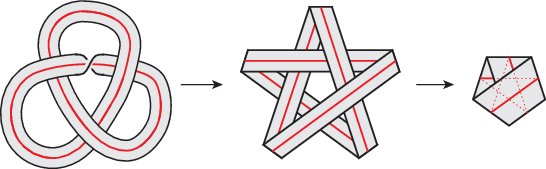}
	\caption{A folded ribbon knot of the trefoil knot}
	\label{fig:folded}
\end{figure}

Kauffman~\cite{K} introduced a {\it folded ribbon knot\/}, which represents a knot by folding a long and thin rectangular strip to become flat.
By the construction of a folded ribbon knot, the core of the ribbon is piecewise linear as drawn in Figure~\ref{fig:folded}.
A ribbonlength can be defined as the minimum length of a long and thin strip of paper required to tie the given knot type.
This definition gives us the intuitive notion of tightness of knots, providing a quantitative measure to classify type of knots.
Kauffman~\cite{K} proposed ribbonlength as drawn in Figure~\ref{fig:folded}.
The following definitions were introduced in~\cite{DHLM,DKTZ}.

\begin{definition}
Let $K$ be the knot. 
Then $K_w$ is a {\it folded ribbon knot\/} of $K$  with width $w$ that satisfies following two conditions:
\begin{enumerate}
\item The ribbon is flat and its fold lines are disjoint.
\item The core of $K_w$ is a union of a finite number of consecutively straight lines with crossing information which represents the knot type of $K$.
\end{enumerate}
\end{definition}

To define the ribbonlength of a knot or link, the infimum is used rather than the minimum.
To see why, consider the trivial knot, which can be represented by folding a ribbon along two parallel lines perpendicular to the core.
Then the length of the such ribbon can be decreased to an arbitrarily small.

\begin{definition}
Take a knot (or link) $K$ and a positive real number $w$.
Then
$$\text{Rib}(K)=\inf_{K' \in [K]_w}\frac{\text{Len}(K'_w)}{w}$$
is called a {\it ribbonlength\/} of $K$
where $[K]_w$ is the set of knots (or links) that are equivalent to $K$ and have a folded ribbon knot $K'_w$.
\end{definition}

The upper bound of a ribbonlength of knots has been actively researched.
Tian~\cite{T} showed that the upper bound of knots is
$$\text{Rib}(K) \leq 2c(K)^2 + 6c(K) + 4.$$
and Denne~\cite{D} improved the result.
She showed that for any knot or link $K$, 
$$\text{Rib}(K) \leq 72 c(K)^{3/2} + 32 c(K) + 12 c(K)^{1/2}+4.$$
Furthermore if $K$ is minimally Hamiltonian, then
$$\text{Rib}(K) \leq 9 c(K)^{3/2} + 8 c(K) + 6 c(K)^{1/2}+4.$$
There are many result for the upper bound of ribbonlength of various families of knots~\cite{DHLM,KNY1,KNY2}

In this paper, we give an improved upper bound on the ribbonlength of knots which is linear in crossing number.

\begin{theorem} \label{thm:main}
For any knot or link $K$,
$$\emph{\text{Rib}}(K) \leq 2.5 c(K)+1.$$
\end{theorem}

Denne~\cite{D} mentioned the {\it ribbonlength crossing number problem\/} which asks us to find constants $c_1$, $c_2$, $\alpha$, $\beta$ such that
$$c_1 \cdot c(K)^{\alpha} \leq \text{Rib}(K_w) \leq c_2 \cdot c(K)^{\beta}.$$
She also mentioned that Diao and Kusner conjectured that $\alpha=\frac{1}{2}$ and $\beta=1$.
The authors~\cite{KNY2} found the partial answer of the conjecture that $\beta=1$ for a 2-bridge knot or link.
Theorem~\ref{thm:main} answers a long held conjecture about the relationship between ribbonlength and crossing number for an arbitrary knot or link.

In Section~\ref{sec:pre}, we introduce a binary grid diagrams of knots and links, and a bisected vertex leveling of a plane graph. 
In Section~\ref{sec:construction}, we explain how to construct the folded ribbon knot, and give lemmas to prove the main theorem.
Finally, in Section~\ref{sec:riblen}, we find an upper bound of the ribbonlength of knots and links.

\section{preliminary}\label{sec:pre}

In this section we define a binary grid diagram.
Furthermore, we explain the basic concept of bisected vertex leveling and introduce its relationship with a binary grid diagram.

\subsection{Binary grid diagrams} \label{subsec:bgd}\ 

A {\it grid diagram\/} is a knot diagram consisting of finite number of vertical line segments and the same number of horizontal line segments such that each vertical line segment is placed over the horizontal line segments.
We introduce a special type of a grid diagram as the following definition.
\begin{definition}
    A {\it binary grid diagram\/} is a grid diagram so that each horizontal line segment crosses at most one vertical line segment.
\end{definition}
For convenience in recognizing the figures during the proof process, we assume that the horizontal line segment of the binary grid diagram passes over the vertical line segment.

\begin{figure}[h!]
	\includegraphics{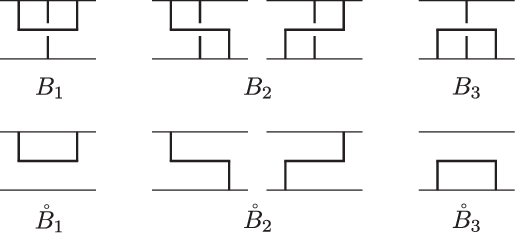}
	\caption{Six types of blocks in a binary grid diagram}
	\label{fig:BGD_block}
\end{figure}

A binary grid diagram can be divided into blocks by cutting horizontally between every consecutive horizontal line segments.
That is, each block contains exactly one horizontal line segment and parts of several vertical line segments.
Thus each block is represented by six types as drawn in Figure~\ref{fig:BGD_block}, according to the shape of vertical line segments which meet the horizontal line segment.
Such blocks are denoted by $B_1$, $B_2$, $B_3$, $\mathring{B}_1$, $\mathring{B}_2$ and $\mathring{B}_3$.

\subsection{Bisected vertex leveling} \label{subsec:vl}\

Let $G$ be a connected plane graph which does not have loops.
A {\it bisected vertex leveling\/} is its ambient isotopy which satisfies the following three conditions;
\begin{enumerate}
    \item Each vertex of $G$ lies between two consecutive horizontal line.
    \item Each edge of $G$ has no maxima and minima as critical points of the height function given by the vertical direction, except its endpoints (vertices).
    \item Each horizontal line cuts G into two pieces, each of which is connected.
\end{enumerate}

The authors~\cite{NOY} first introduced a bisected vertex leveling, and showed the following proposition.
\begin{proposition}[Theorem~1. in \cite{NOY}]\label{prop:bvl}
Let $G$ be a connected plane graph without loops. 
If it has no cut vertex, then $G$ has a bisected vertex leveling.  
\end{proposition}

A knot projection can be interpreted as a 4-valent plane graph.
Let $G$ be a such 4-valent plane graph.
If the given knot projection does not have any nugatory crossing, then $G$ has no loop and no cut vertex.
Therefore, $G$ has a bisected vertex leveling by the following corollary.
\begin{corollary}\label{cor:bvl}
    For any knot diagram with no nugatory crossing, there is a planar isotopy that transforms the knot diagram to a bisected vertex leveling.
\end{corollary}
Note that any minimal crossing diagram of knots has no nugatory crossing.
Therefore, any plane graph corresponding to a knot diagram with minimal crossing has a bisected vertex leveling by Corollary~\ref{cor:bvl}.

\begin{figure}[h!]
	\includegraphics[scale=.8]{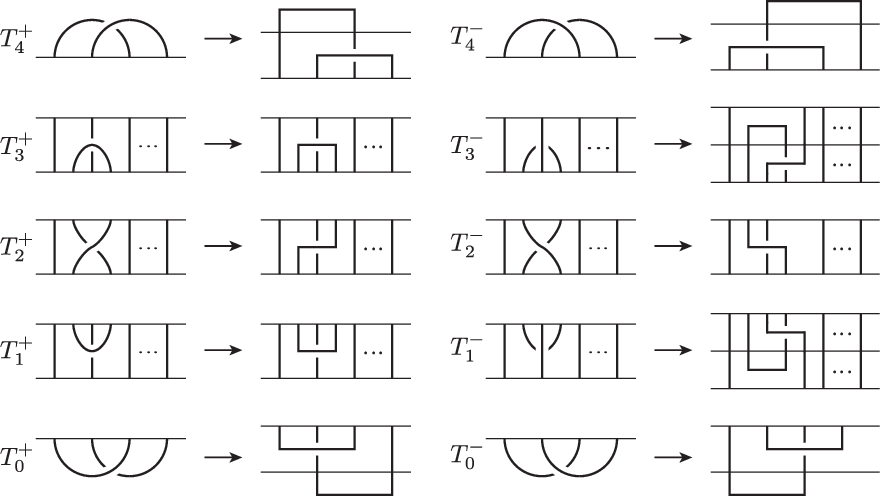}
	\caption{Portions $T_i^{\ast}$'s and their corresponded blocks}
	\label{fig:vl}
\end{figure}

For the convenience of proofs, we label the portion between consecutive horizontal lines for the diagram of recovering crossing information from $G$ as $T_0^+, \dots , T_4^+, T_0^-, \dots , T_4^-$ as drawn in Figure~\ref{fig:vl}.
Moreover, The portion $T_i^+$ and $T_i^-$ are collectively denoted as $T_i^{\ast}$.
By the conditions of a bisected vertex leveling, there are exactly one $T_4^{\ast}$ and exactly one $T_0^{\ast}$.
This implies that the number of $T_3^{\ast}$ is equal to the number of $T_1^{\ast}$.

\begin{remark}
Each portion $T_i^{\ast}$ can be represented by the combination of blocks as the following.
\begin{center}
    \begin{tabular}{lcl}
        $T^{\ast}_0$ & $\mapsto$ & a combination of blocks of types $B_1$ and $\mathring{B}_1$\\
        $T^+_1$ & $\mapsto$ & a block of type $B_1$\\
        $T^-_1$ & $\mapsto$ & a combination of blocks of types $\mathring{B}_1$ and $B_2$\\
        $T^{\ast}_2$ & $\mapsto$ & a block of type $B_2$\\
        $T^+_3$ & $\mapsto$ & a block of type $B_3$\\
        $T^-_3$ & $\mapsto$ & a combination of blocks of types $B_2$ and $\mathring{B}_3$\\
        $T^{\ast}_4$ & $\mapsto$ & a combination of blocks of types $B_3$ and $\mathring{B}_3$\\[2.5pt]
    \end{tabular}
\end{center}
\end{remark}
Therefore any minimal crossing diagram of knots can be represented as a binary grid diagram.

\section{Construction} \label{sec:construction}

In this section, we explain the construction of a folded ribbon knot from a binary grid diagram.
We consider the process of transforming blocks to suitable portions of folded ribbon, and connecting them in the correct order to form the complete folded ribbon knot.
The ribbonlength is determined by specific blocks within a binary grid diagram.
We introduce some lemmas to derive the main theorem.

\begin{lemma} \label{lem:switch}
Let $D$ be a binary grid diagram for knot type of $K$.
If there is a pair of consecutive blocks in which the upper block is of type $B_1$ or $\mathring{B}_1$, and the lower block is of type $\mathring{B}_3$, then we can switch types of the upper and lower blocks while preserving the knot type of $K$.
\end{lemma}

\begin{proof}    
Suppose that there is a pair of consecutive blocks in which the upper block is of type $B_1$ or $\mathring{B}_1$, and the lower block is of type $\mathring{B}_3$.
Let the horizontal line segments in the upper block and the lower block be $[u_1,u_2] \times \{ n+1 \}$ and $[l_1,l_2] \times \{ n \}$, respectively.
If $[l_1,l_2] \cap [u_1,u_2]$ is empty, then switch $y$-coordinates of the upper and the lower horizontal line segments by using the planar isotopy.
Thus, types of the upper and the lower blocks are switched after the above process.

\begin{figure}[h!]
	\includegraphics[scale=.9]{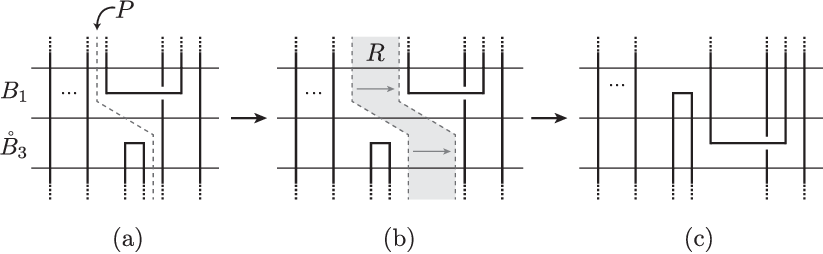}
\caption{A switching process of the upper and lower blocks}
	\label{fig:switch}
\end{figure}

Now we consider the case that $[l_1,l_2] \cap [u_1,u_2]$ is not empty.
First assume that the upper block is of type $B_1$.
Let $c$ be an $x$-coordinate of the vertical line segment which is passing through the horizontal line segment in $B_1$.
Since the lower horizontal line segment does not cross this vertical line segment, there are two cases that $l_2 \in [u_1,c]$ and $l_1 \in [c,u_2]$.

For the first case $l_2 \in [u_1,c]$, we take the piecewise linear line $P$ as drawn in Figure~\ref{fig:switch}(a).
In detail, $P$ consists of two rays $\{ u_1-\frac{1}{2} \} \times [n+1,\infty)$ and  $\{ l_2+\frac{1}{2} \} \times (-\infty,n]$, and the line segment between them.
We partition the plane into three distinct regions, namely $R$, $R_{-}$, and $R_{+}$, based on the $\varepsilon$-tubular neighborhood $R$ of $P$. 
Here, $R_{-}$ and $R_{+}$ correspond to the left and right regions of the complement of $R$, respectively.
Through a planar isotopy, we expand $R$ and translate $R_{+}$ horizontally until the boundary between $R$ and $R_{+}$ is translated by ${l_2-u_1+1}$ in the $x$ direction while keeping $R_{-}$ fixed as drawn in Figure~\ref{fig:switch}(b).
Now we adjust the horizontal spacing of the diagram so that the vertical line segments are spaced by units apart. 
Then we can switch types of the upper and lower blocks similar with the case of $[l_1,l_2] \cap [u_1,u_2] = \phi$ as drawn in Figure~\ref{fig:switch}(c).
For the second case $l_1 \in [c,u_2]$, we take the piecewise linear line $P$ which consists of two rays $\{ u_2+\frac{1}{2} \} \times [n+1,\infty)$ and  $\{ l_1-\frac{1}{2} \} \times (-\infty,n]$, and the line segment between them.
We can also switch the types of the upper and lower blocks similar with the first case.

It remains to consider the case that the upper block is of type $\mathring{B}_1$.
To address this case, we use the same argument in any case where the upper block is of type $B_1$.
Then we can switch the types of the upper and lower blocks while preserving the knot type of $K$.
\end{proof}

\begin{figure}[h!]
	\includegraphics[scale=1]{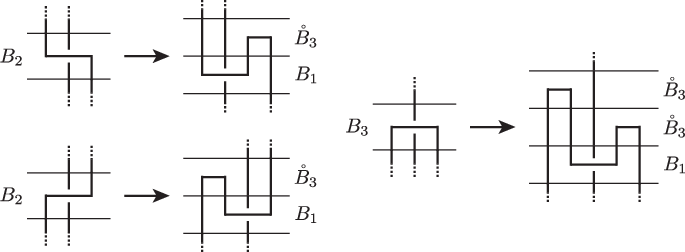}
	\caption{Transforming blocks of types $B_2$ and $B_3$}
	\label{fig:trans}
\end{figure}

We can transform blocks of types $B_2$ and $B_3$ to combinations of blocks of types $B_1$ and $\mathring{B}_3$ as drawn in Figure~\ref{fig:trans}.
So we obtain the following lemma.

\begin{lemma} \label{lem:convert}
Let $D$ be a binary grid diagram for knot type of $K$.
Then there exists an algorithm that which converts $D$ into another binary grid diagram consisting of blocks of types $B_1$, $\mathring{B}_1$, and $\mathring{B}_3$ while maintaining the total number of all blocks of types $B_1$, $B_2$, $B_3$ and $\mathring{B}_1$.
\end{lemma}

\begin{figure}[h!]
	\includegraphics[width=0.8\textwidth]{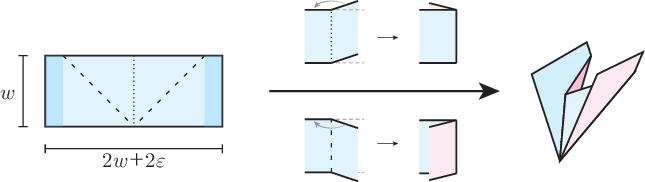}
	\caption{A structure of a paper plane}
	\label{fig:paper}
\end{figure}

Now we introduce a particular object to construct the folded ribbon knot.
A {\it paper plane\/} is a folded ribbon obtained by folding a piece of ribbon three times as shown in Figure~\ref{fig:paper}.
In the figure on the left, the darker regions on the two ends of the paper plane can be made as long or as short as needed in order to join different paper planes together. 
These regions are called {\it wings\/}. 
The horizontal line segment of a block of type $B_1$ or $\mathring{B_1}$ can be represented by a paper plane, while the vertical line segments correspond to wings.

\begin{lemma} \label{lem:count}
Let $D$ be a binary grid diagram for knot type of $K$.
Then
$$\emph{\text{Rib}}(K) \leq 2 ( b_1 + b_2 + b_3 + \mathring{b}_1 ),$$
where $b_i$ and $\mathring{b}_1$ are the number of blocks of type $B_i$ and $\mathring{B}_1$ in $D$, respectively.
\end{lemma}

\begin{figure}[h!]
	\includegraphics[scale=.7]{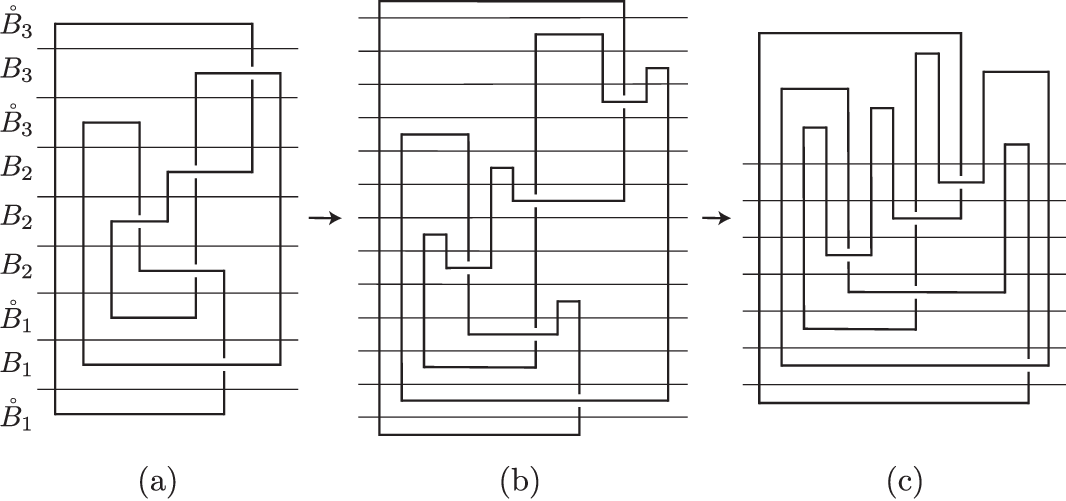}
	\caption{The converted binary grid diagram $D'$}
	\label{fig:bgd}
\end{figure}

\begin{proof}
Let $D$ be a binary grid diagram of a knot $K$ as drawn in Figure~\ref{fig:bgd}(a).
By Lemma~\ref{lem:convert}, We convert all blocks of types $B_2$ and $B_3$ to combinations of blocks of types $B_1$ and $\mathring{B}_3$ as drawn in Figure~\ref{fig:bgd}(b).
Then rearrange blocks such that all blocks of type $\mathring{B}_3$ are lying on the above the other types by Lemma~\ref{lem:switch} as drawn in Figure~\ref{fig:bgd}(c).
We denote $D'$ for the result diagram.
Let the upper part of $D'$ be the part consisting of all blocks of type $\mathring{B}_3$, and the lower part of $D'$ be the other.

Now we represent a folded ribbon knot with width $w$ from $D'$.
We divide the remaining proof into three parts.

\bigskip
\noindent
{\bf Step I.} Transforming the lower part of $D'$ into paper planes

The $i$-th {\it stack\/} is a part of the rearranged $D'$ which consists of $i$ blocks from the bottommost block.
We denote the $i$-th stack as $S_i$.
To construct a ribbon knot, we first transform all blocks in the lower part of $D'$ into paper planes.
We use a sequential process superimposing a new paper plane to the former stack.
So we need the empty space for inserting new paper plane, call this space a {\it insertable space\/}.

We aim to demonstrate, through induction, that the $k$-th stack $S_k$ can be represented by a particular arrangement of $k$ paper planes. 
In this arrangement, all wings are oriented upward, and there are insertable spaces between every pair of consecutive wings for the next paper plane. 
Such an arrangement is called a {\it pile of $k$ paper planes\/}.
The bottommost layer is always composed of $\mathring{B}_1$, and there exists an insertable space between the wings of $\mathring{B}_1$, establishing the initial step.

Now assuming that $S_k$ has been represented by a pile of $k$ paper planes and there is an insertable space between every pair of consecutive wings for the next paper plane.
And then we proceed to implement the appearance up to the $(k+1)$-th block by adding one more paper plane.

Suppose that the $(k+1)$-th block is of type $B_1$. 
Then its horizontal line segment only crosses a single vertical line segment.
There is a wing in the a pile of $k$ paper planes, say the $i$-th wing in order from left to right, which is corresponding to the such vertical line segment.
We insert each wing of the new paper plane into the insertable space.
In detail, the left wing is inserted between the $(i-1)$-th wing and the $i$-th wing, and the right wing is inserted between the $i$-th wing and the $(i+1)$-th wing as drawn in Figure~\ref{fig:ind}.
Here if $i$-th wing is the leftmost(resp. rightmost) wing, then the left(resp. right) wing of the new paper plane is inserted in the left(resp. right) space of the a pile of $k$ paper planes.
Thus we have constructed the pile of $k+1$ paper planes which represents $S_{k+1}$.
We note that this new pile has $2k+2$ wings, and there is an insertable space between each consecutive wing.

\begin{figure}[h!]
	\includegraphics[scale=1]{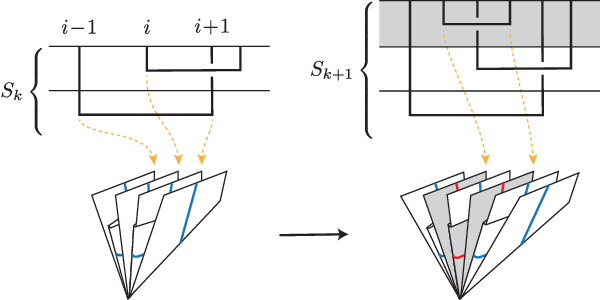}
    \caption{Constructing $S_{k+1}$ from $S_k$ by inserting a paper plane}
	\label{fig:ind}
\end{figure}

It remains to consider that the $(k+1)$-th block is of type $\mathring{B}_1$.
Then there are two cases according to the location of the horizontal line segment of $(k+1)$-th block.
The first case is that the horizontal line segment is lying between two consecutive vertical line segments of $S_k$.
Then there are two wings in the a pile of $k$ paper planes, say the $i$-th and the $(i+1)$-th wing, which are corresponding to the such vertical line segments.
We insert both wings of the new paper plane into the insertable space between the $i$-th wing and the $(i+1)$-th wing.
The second case is that the horizontal line segment are lying in the left(resp. right) of the leftmost(resp. rightmost) vertical line segments of $S_k$.
Then we attach the new paper plane to the left(resp. right) of the a pile of $k$ paper planes.
So we obtain the a pile of $k+1$ paper planes which is represented $S_{k+1}$.
We note that this new pile has $2k+2$ wings, and there is an insertable space between each consecutive wing.

\bigskip
\noindent
{\bf Step II.} Construction of a folded ribbon knot from $D'$

We transform the upper part of $D'$ into the product of the set of arcs in this part with an interval of width $w$ as drawn in Figure~\ref{fig:upper}.
Then it is represented as a union of ribbons with the same number of paper planes. 
We connect the upper part and the lower part from left to right sequentially to construct a folded ribbon knot.
This process guarantees that the core of the resulting folded ribbon knot is same with the diagram $D'$.

\begin{figure}[h!]
	\includegraphics[scale=.93]{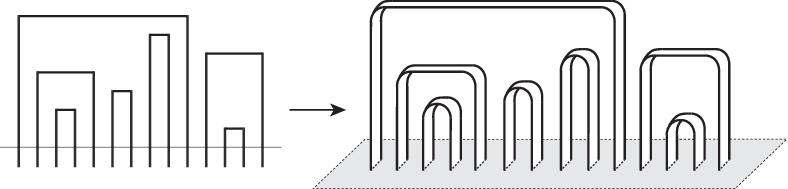}
	\caption{Transforming the upper part of $D'$}
	\label{fig:upper}
\end{figure}

\bigskip
\noindent
{\bf Step III.} $\displaystyle{\text{Rib}(K) \leq 2 (b_1 + b_2 + b_3 + \mathring{b}_1)}$ where $b_i$ and $\mathring{b}_1$ are the number of blocks of type $B_i$ and $\mathring{B}_1$ in $D$, respectively.

Note that both two binary grid diagram $D$ and $D'$ have the same knot type of $K$.
So we check the ribbonlength of a folded ribbon knot which is constructed in Step II.
This folded ribbon knot consists of paper planes representing blocks of types $B_1$ and $\mathring{B}_1$, and their connections.
Since the length of each wing and each connection can be made sufficiently small, the length of the folded ribbon can be reduced to the total length of paper planes.
Note that the length of each paper plane is equal to 2.
Since the folded ribbon knot is constructed by $D'$, it contains the same number of paper planes as the total number of blocks of type $B_1$ and $\mathring{B}_1$ in $D'$.
This also equals the total number of blocks of type $B_1$, $B_2$, $B_3$ and $\mathring{B}_1$ in $D$ by Lemma~\ref{lem:convert}.
Thus a folded ribbon knot of $K$ is represented by $\text{Rib}(K) \leq 2(b_1+b_2+b_3+\mathring{b}_1)$.
\end{proof}

\section{An upper bound of ribbonlength} \label{sec:riblen}

In this section, we find an upper bound on the ribbonlength of knots and links in terms of the crossing number of knots and links, thus proving Theorem~\ref{thm:main}.
Let $K$ be a knot or link.
Then $K$ can be represented as a minimal crossing diagram which includes exactly one $T^{\ast}_0$ and exactly one $T^{\ast}_4$ by using the bisected vertex leveling.
As mentioned in Subsection~\ref{subsec:vl}, all portions $T_i^{\ast}$ are transformed into combinations of blocks to represent $K$ as a binary grid diagram.
We remark that ribbonlength is only affected by the number of blocks of types $B_1$, $B_2$, $B_3$ and $\mathring{B}_1$.
Thus a bisected vertex leveling with smaller number of $T^-_1$ is efficient to find the sharper upper bound of ribbonlength, since it has exactly one $T_0^{\ast}$.

\begin{figure}[h!]
	\includegraphics[scale=1]{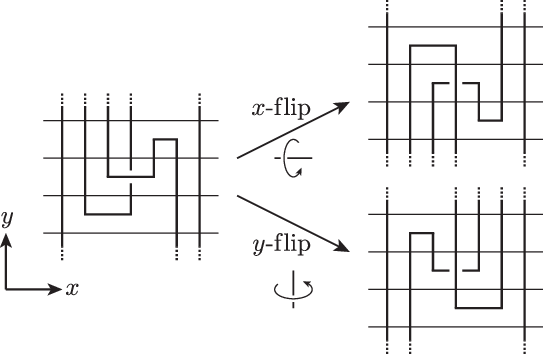}
	\caption{Rotating transformations the $x$-flip and the $y$-flip}
	\label{fig:flip}
\end{figure}

Now we introduce the transformations the $x$-flip and the $y$-flip for a diagram obtained from the bisected vertex leveling of $K$ as drawn in Figure~\ref{fig:flip}.
Here, we consider the diagram is embedded in the three dimensional space $\mathbb{R}^3$, and also consider the rotations in this space.
The $x$-flip and the $y$-flip denote transformations of a knot diagram that rotates about the $x$-axis and the $y$-axis through $\pi$, respectively.
The $x$-flip and the $y$-flip changes the crossing information.
The $x$-flip exchanges $T^{\ast}_i$ and $T^{\ast}_{4-i}$, and the $y$-flip exchanges $T^{\ast}_i$ and $T^{\ast}_{i}$ for $i=1 \dots, 4$.
Especially the $x$-flip and the $y$-flip exchange $T^{\ast}_1$ and $T^{\ast}_3$ as the following table.
\begin{center}
    \begin{tabular}{c|c}
        \ Type of flips \ &   Exchanging operation \\[1.5pt]
\hline
        $x$-flip &  \ \  $T^+_1 \leftrightarrow T^-_3$, \ $T^-_1 \leftrightarrow T^+_3$\ \ \\[2pt]
        $y$-flip &  \ \ $T^+_1 \leftrightarrow T^-_1$, \ $T^+_3 \leftrightarrow T^-_3$\ \ \\[2pt]
    \end{tabular}
\end{center}
This implies that we can change the knot diagram using the $x$-flip and the $y$-flip such that the smallest number among $T^+_1$, $T^-_1$, $T^+_3$ and $T^-_3$ becomes $T^-_1$.
Since the knot diagram includes exactly one each of $T^{\ast}_0$ and $T^{\ast}_4$, the $x$-flip and the $y$-flip do not change each of those numbers.
This implies that the number of $T^-_1$ is at most $\big\lfloor \frac{c(K)-2}{4} \big\rfloor$ in the changed knot diagram.
By the remark in Subsection~\ref{subsec:vl}, each portion of $T^{\ast}_0$ and $T^-_1$ is transformed into two blocks among types $B_1$, $B_2$, $B_3$ and $\mathring{B}_1$.
Moreover each of the other portions is transformed into one or two blocks which contains exactly one block among these four types.
Thus the binary grid diagram obtained from this knot diagram includes at most $c(K)+1+\big\lfloor \frac{c(K)-2}{4} \big\rfloor$ blocks among types $B_1$, $B_2$, $B_3$ and $\mathring{B}_1$.
By Lemma~\ref{lem:count},
\begin{align*}
\text{Rib}(K) & \leq 2(b_1+b_2+b_3+\mathring{b}_1)\\
& \leq 2 \Big( c(K)+1+ \Big\lfloor \frac{c(K)-2}{4} \Big\rfloor \Big)\\
& \leq 2.5c(K)+1,    
\end{align*}
where $b_i$ and $\mathring{b}_1$ are the number of blocks of type $B_i$ and $\mathring{B}_1$ in $D$, respectively.

\bibliography{ribgen.bib} 
\bibliographystyle{abbrv}

\end{document}